\documentclass[12pt]{amsart}
\usepackage{amsfonts}
\usepackage{geometry,amssymb}

\newtheorem{thm}{Theorem}[section]
\newtheorem{prop}[thm]{Proposition}
\newtheorem{cor}[thm]{Corollary}

\theoremstyle{definition}
\newtheorem{defn}[thm]{Definition}

\theoremstyle{remark}

\newcommand{\dist}{{\rm dist}}

\newcommand{\NN}{\mathbb N}
\newcommand{\RR}{\mathbb R}
\newcommand{\CC}{\mathbb C}

\newcommand{\RRd}{{{\mathbb R}^d}}

\newcommand{\LL}{\mathcal L}
\newcommand{\PP}{\mathcal P}
\newcommand{\PPn}{{\mathcal P}_n}
\newcommand{\KK}{\mathcal K}

\newcommand{\bfx}{{x}}

\newcommand{\bfy}{{y}}

\newcommand{\bfa}{{a}}

\newcommand{\bfr}{{\bf r}}

\def\sbt{\subset}
\def\ol{{\bf 0}}
\def\la{\lambda}

\def\pa{\partial}
\def\al{\alpha}

\def\x{{\bf x}}
\def\vp{\varphi}

\def\iy{\infty}

\def\dist{\mathop{\rm dist}}
\def\intt{\mathop{\rm int}}
\def\ol{\overline}
\def\<{\langle}
\def\>{\rangle}
\def\ey{\emptyset}
\def\bcapi{\mathop{\hbox{$\bigcap$}}\limits}

\newcounter{oldresult}
\def\theoldresult{\Alph{oldresult}}
\newenvironment{oldresult}{
  \em
  \vskip 0.10in
  \refstepcounter{oldresult}
  \noindent{\bf Theorem\ \theoldresult.}
}{\vskip 0.10in}

\newcounter{oldprop}
\def\theoldprop{\Alph{oldprop}}

\newcounter{oldlemma}
\def\theoldlemma{\Alph{oldlemma}}
\newenvironment{oldlemma}{
  \em
  \vskip 0.10in
  \refstepcounter{oldlemma}
  \noindent{\bf Lemma\ \theoldlemma.}
}{\vskip 0.10in}

\newcounter{oldcor}
\def\theoldcor{\Alph{oldcor}}

\newcounter{oldconjecture}
\def\theoldconjecture{\Alph{oldconjecture}}
\newenvironment{oldconjecture}{
  \em
  \vskip 0.10in
  \refstepcounter{oldconjecture}
  \noindent{\bf Conjecture\ \theoldconjecture.}
}{\vskip 0.10in}

\newcounter{hypothesis}
\def\thehypothesis{\Alph{hypothesis}}
\newenvironment{hypothesis}{
  \em
  \vskip 0.10in
  \refstepcounter{hypothesis}
  \noindent{\bf Hypothesis\ \thehypothesis.}
}{\vskip 0.10in}

\begin{document}
\title[Polynomial inequalities on normed spaces] {Some polynomial
inequalities on real normed spaces}

\author{ Szil{\'a}rd\ Gy. R{\'e}v{\'e}sz}





\begin{abstract} We consider various inequalities for polynomials,
with an emphasis on the most fundamental inequalities of
approximation theory. In the sequel a key role is played by the
generalized Minkowski functional $\alpha(K,x)$, already being used
by Minkowski and contemporaries and having occurred in
approximation theory in the work of Rivlin and Shapiro in the
early sixties. We try to compare real, geometric methods and
complex, pluripotential theoretical approaches, where possible,
and formulate a number of questions to be decided in the future.
An extensive bibliography is given to direct the reader even in
topics we do not have space to cover in more detail.
\end{abstract}

\maketitle

\section{Introduction}\label{sec:intro}

\subsection{}\label{sec:Chebypoly} In our present work, as well as throughout
and all over approximation theory, a distinguished role is played
by the (univariate) Chebyshev polynomials of the first kind. These
can be defined as
\begin{equation}\label{Tndef}
\begin{array}{rl} T_n(x)&:= 2^{n-1}\prod\limits_{j=1}^n
\left(x-\cos\left({\frac{(2j-1)\pi}{2n}}\right)\right)=
\\[4mm]
&= {\frac{1}{2}}\left\{{(x+\sqrt{x^2-1})}^n
+{(x-\sqrt{x^2-1})}^n\right\},
\end{array}
\end{equation}
while the most used expression for them is the first part of the
formula
\begin{equation}\label{Tncos}
T_n(x) = \begin{cases} \cos(n\arccos x), & \vert x \vert \le 1\\
{\rm{sgn}}(x)^n \cosh(n\cosh^{-1} \vert x\vert ), &\vert x\vert
\ge 1 ~.
\end{cases}
\end{equation}
It is well-known that regarding the modulus of
$p(x)/\|p\|_{[-1,1]}$, $T_n$ is extremal simultaneously for all
$x$ with $|x|>1$. For this we refer to \cite[Theorem 1.2.2,
Chapter 5]{MMR} or \cite[(2.37), p.~108]{Ri}. The very same
Chebyshev polynomial has many extremal properties, see for
instance \cite[\S 2.7]{Ri}. In particular, it is also extremal
concerning its ``speed of growth towards infinity'', which can be
precisely described by its leading coefficient. Since for all
$p(x)= \sum_{k=1}^n a_kx^k \in\PPn(\RR)$ we have
$a_n=\lim_{x\to+\infty} {p(x)}/{x^n}$, a polynomial is extremal
concerning its growth towards infinity iff the leading coefficient
is extremal. Thus, for
\begin{equation}\label{univleadingc}
\max \left\{ a_n: p(x)= \sum_{k=0}^n a_kx^k \in\PPn(\RR),\ \Vert
p\Vert_{[-1,1]}\le 1\right\},
\end{equation}
Chebyshev's polynomial is again the extremal case.

\subsection{}\label{sec:poly}
In the present survey we focus on extensions of the most
well-known and classical inequalities of approximation theory for
algebraic polynomials on $\RR$ to the case of infinitely many
variables, i.e. to normed spaces. In all what follows, $X$ is a
real normed space, $X^*=\LL(X,\RR)$ is the usual dual space, and
$S:=S_X$, $S^*:=S_{X^*}$, $B:=B_X$ and $B^*:=B_{X^*}$ are the unit
spheres and (closed) unit balls of $X$ and $X^*$, respectively.
Moreover, $\PP=\PP(X)$ and $\PPn=\PPn(X)$ will denote the space of
continuous (i.e., bounded) polynomials of free degree and of
degree at most $n$, respectively, from $X$ to $\RR$.

There are several ways to introduce continuous polynomials over
$X$, one being the linear algebraic way of writing
\begin{equation}\label{Pnassum}
\PPn:=\PP^*_0+\PP^*_1+\dots+\PP^*_n ~, \qquad \text{and} \qquad
\PP:=\bigcup_{n=0}^{\infty}\PPn
\end{equation}
with $\PP^*_k$ (or, in another notation, $\PP(^kX;\RR)$) denoting
the space of homogeneous (continuous) polynomials of degree
(exactly) $k\in\NN$. That is, one considers bounded $k$-linear
forms
\begin{equation}\label{Lk}
L\in \LL(X^k\to\RR)
\end{equation}
together with their ``diagonal functions"
\begin{equation}\label{Lhat}
\widehat{L}:X\to\RR,\qquad \widehat{L}(x):=L(x,x,\dots,x)
\end{equation}
and defines $\PP^*_k$ as the set of all $\widehat{L}$ for $L$
running $\LL(^kX):=\LL(X^k\to\RR)$. In fact, it is sufficient to
identify equivalent linear forms (that is, those having identical
diagonal functions) by selecting the unique symmetric one among
them: in other words, to let $L$ run over $\LL^s(^kX)$ denoting
(real) {\em symmetric} $k$-linear forms. Building up the notion of
polynomials that way is equivalent to the definition
\begin{equation}\label{PnX}
\PPn:=\left\{p:X\to\RR~:~\|p\|<\infty, ~p|_{Y+y}\in\PPn(\RR)
~\text{for all}~~ Y\le X,~\dim Y =1,~ y\in X \right\}
\end{equation}
or to the definition arising from combining \eqref{Pnassum} and
\begin{equation}\label{PkXstar}
\PP^*_k:=\left\{p:X\to\RR~:~\|p\|<\infty, ~p|_{Y} \in\PPn^*(\RR^2)
~\text{for all}~~ Y\le X,~\dim Y =2 \right\}~.
\end{equation}
Here and throughout the paper for any set $K\subset X$ and
function $f:X\to\RR$ we denote, as usual,
$$
\|f\|_K:=\sup_K |f| \qquad \text{and}\qquad \|f\|:=\|f\|_B~.
$$
For equivalent definitions of and introduction to polynomials over
real normed spaces see \cite[Chapter 1]{Din} and also \cite{Coat,
H, Kell, T}. In particular, it is well-known that
\begin{equation}\label{LhatLLhat}
\|\widehat{L}\|\le\|L\|\le C(n,X) \|\widehat{L}\|~\qquad
{\text{for all}}\qquad L\in\LL^s(^n X)~,
\end{equation}
and that $C(n,X)\le n^n/n!$, \cite{Din}, while $C(n,X)=1$ if $X$
is a Hilbert space (Banach's Theorem, see \cite{Banach, Din, H0}).
Similarly to \eqref{LhatLLhat}, one can consider special
homogeneous polynomials which can be written as products of linear
forms, i.e. $L(x_1,x_2,\dots,x_n)=\prod_{j=1}^n f_j(x_j)$ with
$f_j\in\LL(X\to \RR)$. Then $\|L\|=\prod_{j=1}^n\|f_j\|$, i.e. the
product of the norms, and one compares to the norm of the
corresponding homogeneous polynomial, i.e. to
$\|\widehat{L}\|=\|\prod_{j=1}^n f_j\|$. Note that here $L$ is far
from being symmetric, and this yields to an essentially different
question, with the similarly defined polarization constants now
ranging up to $n^n$, see e.g \cite{BALL1, BST, Ry}.

These polarization problems are typical examples of genuinely
multivariate inequalities, as in dimension 1 they simply
degenerate. Since our focus is different, we direct the reader's
attention to \cite{ARIAS, BALL1, BALL2, BALL3, BST, GM, GV, Mate2,
PR, KRS} and also to \cite{Mate} right in this volume. However,
passing by we note that \eqref{LhatLLhat} already implies that a
(symmetric) $n$-linear form $L$ is bounded iff its diagonal
function -- i.e. the associated homogeneous polynomial defined by
$\widehat{L}$ -- is bounded, and that a polynomial
\begin{equation}\label{passumpstar}
p=p^*_0+p^*_1+\dots +p^*_n \qquad \text{with} \qquad
p^*_k=\widehat{L_k}~,
\end{equation}
$L_k: X^k\to \RR$ being a (symmetric) $k$-linear mapping, is
bounded iff $L_k$ are such for all $k=1,\dots,n$. Hence in all
what follows we are free to talk about boundedness or continuity
of these polynomials without specifying in detail whether $p$, or
$p_k^*$ are assumed to be continuous or bounded.

\subsection{}\label{convexbody} In the following classical
inequalities of approximation theory the usual condition of
normalization is that $\|p\|_I\le 1$, where $p\in\PPn(\RR)$ and
$I=[-1,1]$ (or, sometimes, some other interval $[a,b]$). In $\RR$
all the convex bodies are just intervals, and linear substitution
allows to restrict ourselves to $I$, but in higher dimensions
there is a great variety of convex bodies to deal with. Recall
that a set $K\subset X$ is called {\it convex body} in a normed
space (or in a topological vector space) $X$ if it is a bounded,
closed convex set that has a non-empty interior.

The convex body $K$ is {\em symmetric}, iff there exists a center
of symmetry $x$ so that reflection of $K$ at $x$ leaves the set
invariant, that is, $K=-(K-x)+x=-K+2x$. In the following we will
term $K$ to be {\em centrally symmetric} if it is symmetric with
respect to the origin, i.e. if $K=-K$. This occurs iff $K$ can be
considered the unit ball with respect to a norm $\|\cdot\|_{(K)}$,
which is then equivalent to the original norm $\|\cdot\|$ of the
space $X$ in view of $B_{ X,\, \|\cdot\|}( 0,r)\sbt K\sbt
B_{X,\,\|\cdot\|}( 0,R)$.

The {\it central symmetrization\/} or {\it half difference body\/}
(cf. \cite{HC}, p.~135 and ~362, respectively) of a set $K$ in a
normed space $X$ is
\begin{equation}\label{(2.1)}
C:= C(K):= {\frac{1}{2}}(K-K):= \left\{{\frac{1}{2}}( x- y)\colon
x, y\in K \right\} .
\end{equation}
The central symmetrization of $K$ is centrally symmetric with
respect to the origin. In case $K$ is a convex body, we also have
$ 0\in \intt C$\footnote{Throughout the paper we denote when
convenient $C(K),\, \tau (K,v),\, \al (K, x ), \,w(K,v^*)$, etc.
by $C , \, \tau ,\, \al , \, w$, etc., respectively.}. On the
other hand, even though $K$ is assumed to be closed, $C$ is not
necessarily closed (c.f. \cite[Section 6]{RS}), hence $C$ is not a
convex body in general. Nevertheless, the closure $\ol C$ of $C$
is a symmetric convex body, which is also {\it fat}, and $\intt
C\sbt C\sbt \ol{\intt C}= \ol C$.

The {\it ``maximal chord"} of $K$ in direction of $v\ne  0$ is
\begin{equation}\label{maxchord}
\begin{gathered}
\tau (K,v):= \sup\{ \la \ge 0:\exists \,\,  y,\, z\in K \,\hbox{
s.t. } z= y+\lambda v\}=
\\
=\sup \big\{ \la \ge 0\colon K\cap(K+\la v)\ne \ey\big\}=
\\
=\sup\{\la \ge 0\colon \la v\in K-K\}= 2\sup\big\{ \la>0\colon
\la v \in C\big\}=
\\
=2\max\big\{ \la \ge
0\colon \la v \in \ol{C}\big\}= \tau\big( C, v\big).
\end{gathered}
\end{equation}
Usually $\tau(K, v)$ is not a ``maximal" chord length, but only a
supremum, however we shall use the familiar finite dimensional
terminology (see for example \cite{W}).

The {\it support function\/} to $K$, where $K$ can be an arbitrary
set, is defined for all $ v^*\in X^*$ (sometimes only for $ v^*\in
S^*$) as
\begin{equation}\label{(2.3)}
h(K, v^*):= \sup_K v^*=\sup\big\{\< v^*, x\>\colon  x\in K\big\},
\end{equation}
and the {\it width of $K$ in direction\/} $ v^*\in X^*$ (or $
v^*\in S^*$) is

\begin{equation}\label{(2.4)}
\begin{gathered}
\hfill w(K, v^*):= h(K, v^*)+h(K,- v^*)=\sup_K v^*+\sup_K(- v^*)=
\\
=\sup\big\{\< v^*, x- y\>\colon  x, y\in K\big\}=2h\big(C,
v^*\big)=w\big(C, v^*\big).
\end{gathered}
\end{equation}
Let us introduce the notations
\begin{equation}\label{(2.5)}
X_t( v^*):= \big\{ x\in X \colon \< v^*, x\>\le t\big\}, \qquad
X(K, v^*):=X_{h(K, v^*)}( v^*).
\end{equation}
Clearly the closed halfspace $X(K, v^*)$ contains $K$, and the
hyperplane
\begin{equation}\label{(2.6)}
H(K, v^*):=H_{h(K, v^*)}( v^*), \quad H_t( v^*):= \big\{ x\in X
\colon \< v^*, x\>= t\big\}=\pa X_t( v^*)
\end{equation}
is a supporting hyperplane\footnote{Note that throughout the paper
we mean ``supporting" in the weak sense, that is, we do not
require $K\cap H(K, v^*)\ne\ey$, but only $\dist \big(K,H(K,
v^*)\big)=0$. The same convention is in effect for other
supporting objects as halfspaces, layers etc.} to $K$.

A {\it layer\/} (sometimes also called {\it slab\/}, {\it plank\/}
or {\it strip\/}) is the region of $X$ enclosed by two parallel
hyperplanes, i.e.
\begin{equation}\label{(2.7)}
L_{r,s}( v^*):= \big\{ x\in X \colon r\le \< v^*, x\>\le s\big\}
=X_s( v^*)\cap X_{-r}(- v^*),
\end{equation}
while the {\it supporting layer\/} or {\it fitting layer\/} of $K$
with normal $ v^*$ is
\begin{equation}\label{(2.8)}
\begin{gathered}
L(K, v^*):= X(K, v^*)\cap X(K,- v^*)=L_{-h(K,- v^*),\,h(K, v^*)}(
v^*)=
\\
=\big\{ x\in X \colon -h(K,- v^*)\le \< v^*, x\>\le h(K,
v^*)\big\}.
\end{gathered}
\end{equation}

\subsection{}\label{genMinkowski} In $\RR$ the position of a point
$x\in\RR$ with respect to the "convex body" $I$ can be expressed
simply by $|x|$ (as $\pm x$ occupy symmetric positions). However,
to quantify the position of $x\in X$ with respect to the convex
body $K\subset X$ is a problem of several possible answers. In
this regard the most frequent tool is the Minkowski functional.
For any $ x\in X$ the {\it Minkowski functional\/} or {\it
(Minkowski) distance function\/} \cite[p.~57]{HC} or {\it gauge\/}
\cite[p.~28]{Ro} or {\it Minkowski gauge functional\/} \cite[\S
1.1(d)]{P} is defined as
\begin{equation}\label{(1.1)}
\vp_K( x):= \inf\{\la> 0\colon  x\in\la K\}~.
\end{equation}
Clearly \eqref{(1.1)} is a norm on $ X $ if and only if the convex
body $K$ is centrally symmetric with respect to the origin. If
$K\sbt X $ is a centrally symmetric convex body, then the norm
${\|\cdot\|}_{(K)}:=\vp_K$ can be used successfully in
approximation theoretic questions as well. As said above, for
${\|\cdot\|}_{(K)}$ the unit ball of $X$ will be $K$ itself,
$B_{X,\,{\|\cdot\|}_{(K)}}( 0,1)= K $.

In case $K$ is nonsymmetric, \eqref{(1.1)} still can be used. But
then even the choice of the homothetic centre is questionable
since the use of any alternative gauge functional $$ \vp_{K,
{x_0}}( x):= \inf\big\{\la> 0\colon  x\in \la(K- x_0)\big\} $$ is
equally well justified. Moreover, neither is good enough for the
applications.

One of the key points of these notes is to highlight the role of
the so called {\em generalized Minkowski functional} in the above
quantification problem. It seems that the most appropriate means
to apply in the inequalities of our interest are provided by this
notion. This generalized Minkowski functional $\alpha(K,x)$ also
goes back to Minkowski \cite{Mi} and Radon \cite{Rad}, see also
\cite{Gr}, \cite{RS}. There are several ways to introduce it, but
perhaps the most appealing is the following construction.

By convexity, $K$ is the intersection of its ``supporting
halfpaces" $X(K, v^*)$, and grouping opposite normals we get
\begin{equation}\label{(2.9)}
K=\bcapi_{ v^*\in S^*}X(K, v^*)=\bcapi_{ v^*\in S^*} L(K, v^*).
\end{equation}
Any layer \eqref{(2.7)} can be homothetically dilated with
quotient $\la \ge 0$ at any of its symmetry centers lying on the
symmetry hyperplane $H_{r+\frac{s}{2}}( v^*)$ to obtain
\begin{equation}\label{(2.10)}
L^\la_{r,s}( v^*):= \left\{ x\in X \colon
{\frac{\la+1}{2}}r-{\frac{\la-1}{2}}s\le \< v^*, x\>\le
{\frac{\la+1}{2}}s-{\frac{\la-1}{2}}r\right\}.
\end{equation}
In particular, we have also defined
\begin{equation}\label{(2.11)}
L^\la(K, v^*)=L_{-h(K,- v^*),\,h(K, v^*)}^\la ( v^*)
\end{equation}
and by using \eqref{(2.11)} one can even define
\begin{equation}\label{(2.12)}
K^\la:= \bcapi_{ v^*\in S^*} L^\la(K, v^*) .
\end{equation}
Note that $\;K^\la$ can be empty for small values of $\la$.
Although not needed here, it is worth mentioning that a nice
formula, due to R. Schneider and E. Makai (for $\lambda \ge 1$ and
for $\lambda \le 1$, respectively) states that
\begin{equation}\label{Klamnda}
K^\lambda = \begin{cases} \ol{K+(\lambda -1)C(K)}=\ol{
\frac{\lambda +1}{2} K- \frac{\lambda -1}{2}K}
\\
K \sim (1-\la)\ol{C}=\frac{1+\la}{2}K \sim \frac{1-\la}{2}(-K),
\end{cases}
\end{equation}
see \cite[Propositions 7.1 and 7.3]{RS}, with $\sim$ denoting
Minkowski difference: $ A\sim B:=\{ \x\in X : \x+B\subset A\}$.
The sets \eqref{Klamnda} were first extensively studied by Hammer
\cite{Ha}.

Using the convex, closed, bounded, increasing and (as easily seen,
c.f. \cite[Proposition 3.3 ]{RS}) even absorbing set system
${\{K^\la\}}_{\la\ge 0}$, the {\it generalized Minkowski
functional\/} or {\it gauge functional\/} is defined as
\begin{equation}\label{alphadef}
\al(K, x):= \inf\{\la\ge 0\colon x\in K^\la\}.
\end{equation}

There are other possibilities to define $\alpha(K,x)$
equivalently. First let
\begin{equation}\label{gammadef}
\gamma(K, x) :=  \inf \left\{ 2 \frac{\sqrt{||x-a||~||x-b||}}
{||a-b||}\,: a, b \in \partial{K}, {\rm such}\, {\rm that}~~ x \in
[a,b] \right\}.
\end{equation}
Then we have
\begin{equation}\label{}
\alpha(K, x) = \sqrt{ 1 -\gamma^{2}(K, {x})}.
\end{equation}

Also, one can consider the original definition of Minkowski
\cite{Mi}, as presented in Gr{\"u}nbaum's article,  \cite[p.
246]{Gr}. Denote
\begin{equation}\label{lineart}
t:= t(K, v^*, x ):= \frac{2\< v^*, x \>-h(K, v^*)+h(K,- v^*)}{
w(K, v^*)}.
\end{equation}

For fixed $ v^*$ this function is an affine linear functional in $
x \in X$, while for fixed $ x $ it is a norm-continuous mapping
from $S^*$ (or $X^*\setminus \{ 0 \}$) to $\RR$. In fact, for
fixed $ v^*\in S^*$, $t$ maps the layer $L(K, v^*)$ to $[-1,1]$,
and $L^\eta(K, v^*)$ to $[-\eta,\eta]$. Therefore, the two forms
of the following definition are really equivalent;

\begin{equation}\label{lambdadef}
\begin{aligned}
\la:&= \la(K, x ):= \sup \big\{\eta>0\colon \exists  v^*\in S^*,\;
x \in \pa L^\eta(K, v^*)\big\}
\\
&=\sup\Big\{\big|t(K, v^*, x )\big|\colon  v^*\in S^*\Big\}
=\sup\Big\{ t(K, v^*, x ) \colon  v^*\in S^*\Big\} .
\end{aligned}
\end{equation}
Note that $t(K, v^*, x )=-t(K,- v^*, x )$ and therefore we don't
have to use the absolute value.

In fact, $\la(K, x )$ expresses the supremum of the ratios of the
distances between the point $ x $ and the symmetry hyperplane
$\frac12 (H(K, v^*)+H(K,- v^*))$ of any layer $L(K, v^*)$ and the
half-width $w(K, v^*)/2$. Now Minkowski's definition was $$
\varphi(K, x ):= \inf \left\{ \frac{\min\{\dist\left( x ,H(K,
v^*)\right), \dist\left( x , H(K,- v^*)\right)\}}
{\max\{\dist\left( x , H(K, v^*)\right), \dist\left( x ,H(K,-
v^*)\right)\}} \, : \,  v^*\in S^* \right\}, $$ which clearly
implies the relation $$ \varphi(K, x )=\frac{1-\la(K, x
)}{1+\la(K, x )} \qquad ( x  \in K). $$ Although this $\varphi(K,
x )$ seems to be used traditionally only for $ x \in K$, extending
the definition to arbitrary $ x  \in X$ yields the similar
relation $$ \varphi(K, x )=\frac{|1-\la(K, x )|}{1+\la(K, x )}
\qquad ( x  \in X). $$

Now the above definitions are connected simply by
\begin{equation}\label{alphalambda}
\lambda(K,x)=\alpha(K,x)~.
\end{equation}

In fact, usefulness of \eqref{alphadef} and the possibility of the
wide ranging applications stems from the fact that this geometric
quantity incorporates quite nicely the geometric aspects of the
configuration of $x$ with respect to $K$, which is mirrored by
about a dozen (!), sometimes strikingly different-looking,
equivalent formulations of it. For the above and many other
equivalent formulations with full proofs, further geometric
properties and some notes on the applications in approximation
theory see \cite{RS} and the references therein.

\section{Chebyshev type problems of polynomial growth}
\label{sec:Cheby}

Chebyshev problems are, in fact, a large class of problems. We
select from these only Chebyshev-type extremal problems concerning
growth of real polynomials. There are further questions we do not
address here, one important class being the problem of
approximating a prescribed "main term", i.e. some homogeneous term
or polynomial of given degree $n$, by the collection of lower
degree or lower rank (in lexigographical order) terms. To these
questions we refer to \cite{BC, BHN, Bos, Gear, NX, Peh, Rei1,
Rei2}, and the references therein.

The general question we will be dealing with can be formulated as
follows: ``How large can a polynomial be at a point $ x \in X$, or
when $ x \to\iy$?" More precisely, we are interested in
determining for arbitrary fixed $x \in X$
\begin{equation}\label{Chebyfindef}
C_n(K, x ):= \sup\big\{p( x ):p\in \PPn,\ {\|p\|}_K\le 1\big\},
\end{equation}
or for some $ v\in X $ (say, with $\| v\|=1$)
\begin{equation}\label{Chebyinfindef}
A_n(K, v):= \sup\big\{p_n^*( v)\colon p\in \PPn \hbox{ satisfying
\eqref{passumpstar},~}~ {\|p\|}_K\le 1\big\}.
\end{equation}
Clearly $C_n$ specifies the possible size of a polynomial at a
given point, while $A_n$ is its order of growth towards infinity
in a given direction. Note the appearance of the $n$-homogeneous
part $p_n^*$ in \eqref{Chebyinfindef}. Hence it is apparent that
$A_n$ is a kind of a formulation of the limiting case of $C_n$.
Indeed, it is easy to see by lower order homogeneity of all the
other terms, that for $p$ represented as in \eqref{passumpstar} we
have
\begin{equation}\label{homolimit}
p_n^*(v)= \lim\limits_{\lambda\to+\infty} \frac{p(\lambda
v)}{\lambda^n}~,
\end{equation}
hence a precise knowledge of $p(\lambda v)$ suffices. Both
problems are classical and fundamental in the theory of
approximation, see e.g. \cite{MMR} or \cite{Ri} for the one and a
half century old single variable result and its many consequences,
variations and extensions.

\subsection{}\label{RivlinShapiro}

As we have mentioned in the introduction, even the above
formulation \eqref{(2.12)} and \eqref{alphadef} of the definition
of $\al(K, x )$ was applied first in work on these questions,
particularly on \eqref{Chebyfindef}, where a quantification of the
position of $ x $ with respect to $K$ is needed. To the best of
our knowledge, application of the generalized Minkowski functional
penetrated into approximation theory and polynomial inequalities
first in the fundamental work \cite{RiSh} by Rivlin--Shapiro.
There they proved the following.
\begin{oldresult}{\bf (Rivlin-Shapiro, 1961).}\label{th:RivlinShapiro}
Let $K\subset\RRd$ be a strictly convex body and
$x\in\RRd\setminus K$. Then we have
\begin{equation}\label{CnTn}
C_n(K, x )=T_n\big(\al(K, x )\big).
\end{equation}
Moreover, $C_n(K, x )$ is actually a maximum, attained by
\begin{equation}\label{PTnt}
P( x ):= T_n\big(t(K, v^*, x )\big).
\end{equation}
Here $T_n$ is the classical Chebyshev polynomial \eqref{Tndef},
while $t(K, v^*, x )$ is the linear expression defined in
\eqref{lineart} and $ v^*$ is some appropriately chosen linear
functional from $S^*$.
\end{oldresult}

Note that actually the restriction $ x \notin K$ is natural, as
$p\equiv 1\in \PPn$, and thus for $ x \in K$ we always have $C_n
(K, x )=1$.

Apart from involving the generalized Minkowski functional, Rivlin
and Shapiro naturally used the following helpful auxiliary
proposition from the geometry of $\RRd$.

\begin{oldlemma}{\bf (Parallel supproting hyperplanes lemma).}
\label{l:parsup}
Let $K\subset\RRd$ be a convex body, and $x\in K$
arbitrary. Then there exists at least one straight line $\ell$
through the point $x$ so that $K\cap\ell=[a,b]$ with some $a\ne b$
and $a,b\in \partial K$ and $K$ has parallel supporting
hyperplanes at $a$ and $b$.
\end{oldlemma}
This standard fact was well-known to geometers for long, and many
authors used it without reference or proof, see eg. \cite[p. 990]
{BANG}, \cite{Ha}, nevertheless, in approximation theory some
reproving occurred later on. It is useful both in proving the
result and to find the extremal polynomial exhibiting exactness of
the upper estimate.

Rivlin and Shapiro assumed strict convexity -- which means that no
straight line segment can lie on the boundary $\partial K$ -- for
they needed it in order to apply their basic method, that of the
extremal signatures. In fact they needed this condition in
proving Lemma \ref{l:parsup} {\em by use of extremal signatures},
which was their goal in illustrating the diverse applications of
the method. Indeed, their method proved to be very successful in
multivariate polynomial problems, but they themselves remarked in
\cite{RiSh} that regarding the Chebyshev problem, a direct, more
geometrical argument can give more\footnote{In fact, they present
this as Problem 3 on pages 694-696 of the paper, and start by
explicitly writing "... this problem ... may also be solved
without the methods of this paper." Then after proving the
assertion of Lemma \ref{l:parsup}, they remark once again: "It
is, of course, possible to obtain this result more geometrically
but with (42) at hand we prefer to utilize it." And ending the
application to Problem 3, they state once again: "To sum up: To
solve Problem 3 we need only a pair of parallel supporting
hyperplanes to $K$ such that the points of tangency, $P_1$ and
$P_2$, are collinear with $P_0$." ($P_0$ in their notation
corresponds to the point $x$ in ours.) Then they describe once
again how the corresponding extremal value and the extremal ridge
Chebyshev polynomial is found once these hyperplanes are given.}.
The -- from here quite straightforward -- proof for the case of a
not necessarily strictly convex $K\subset\RRd$ was then presented
in \cite{KS}.

\subsection{}\label{convbodyapprox}

However, there is no need for any new proof until we keep working
in $\RRd$, as Theorem \ref{th:RivlinShapiro} of Rivlin and Shapiro
for strictly convex bodies directly implies the general case
once we take into account the next standard fact.
\begin{oldlemma}\label{l:strictapprox}{\bf (Convex bodies
approximation lemma).} Any convex body $K\subset\RRd$ can be
approximated arbitrarily closely by strictly convex bodies of
$\RRd$.
\end{oldlemma}
Here, naturally, the approximation is meant in the Hausdorff
distance sense, that is in
\begin{equation}\label{Hausdorff}
\delta(K,M):= \max\left\{ \sup_{x\in K} \inf_{y\in M} \|x-y\|, \
\sup_{y\in M} \inf_{x\in K}\|x-y\|\right\}~.
\end{equation}
This can be a kind of folklore among geometers, but to page out a
proof was difficult. Nevertheless, several ideas of proofs were
suggested by colleagues working in geometry, so it can certainly
not be considered an unknown fact. In fact, e.g. in multivariate
complex analysis this is used even with the stronger requirement
that the approximating convex bodies monotonically decrease to $K$
and have even real analytic boundaries, see e.g. \cite[Proposition
2.2]{BLM}. Anyway, we sketch two proofs in the sequel.

{\em First proof of Lemma \ref{l:strictapprox}.} \, Working in
$\RRd$ one may fix any positive $\epsilon$, approximate the given
convex body $K$ within $\epsilon/3$ by some polyhedra $P$, then
$P$ by a special polyhedra $S$ with only 1-codimensional
simplices as sides, and then finally change very slightly the
(then finitely many) halfspaces, giving as their intersection
$S$, so that the resulting body be strictly convex. A way for
this last change is to substitute each halfspace $Q$ by a large
ball $B$, exhibiting a very small Hausdorff distance
$\delta(Q\cap B(0,R),B\cap B(0,R))$, where $R$ is taken so large
(but fixed) that a given large neighbourhood of $K$ is already
contained in it. In fact, it is easier to see that we get what we
want if we construct these balls the following way. We pick up
one point from the relative interior of each of the sides of $S$,
and move it slightly outward in normal direction: then the
corresponding balls $B$ are defined as the balls drawn around
those simplices of full dimension, which arise from the original
sides and the corresponding, slightly moved points outside.
Clearly if the points to be moved are fixed, and the length of
the move is fixed for all sides equally as, say, $\delta$, then
in function of $\delta\to 0$, the intersection of these balls,
(which always contain $S$ for small enough $\delta$), will finally
shrink to $S$. That concludes the proof of the
lemma\footnote{This nice constructive proof was communicated to
us by Bal{\'a}zs Csik{\'o}s.}.

Having this approximation lemma the proof of the general case of
Theorem \ref{th:RivlinShapiro} is done by referring to the
continuity of $T_n$ and also of $\alpha$, the latter understood as
a function on $\KK \times X$, where $\KK$ denotes the set of all
convex bodies, and is equipped with the metric of Hausdorff
distance \eqref{Hausdorff}. Even the extremal polynomial
\eqref{PTnt} obtains using the corresponding extremal polynomials
of the strict convex case and compactness in $\RRd$.

But is continuity clear? Well, continuity of $\alpha$ can be
checked explicitly, but it may be rather tedious, compared to our
expectations that it should be such anyway. So the best is to get
around any tedious calculations, and prove something even better,
that of {\em convexity} in $x$ and admitting a Lipschitz bound
even as a two-variable function on $\KK \times X $, see
\cite[Theorem 5.5]{RS}. (Actually, here \cite[Lemma 5.4]{RS}
suffices.)

In fact, continuity of $\alpha$ holds even in the normed space
setting, which is something we could not get through compactness
or direct calculations, but by combining convexity, Lipschitz
bounds, and, in view of infinite dimension, even the fact that
$\alpha$ is bounded on bounded sets, which is also necessary, see
\cite{P, Ro}. For the whole assertion see \cite[Corollary
6.1]{RS}.

\subsection{}\label{infinitedimext}

All that raise the question whether we can go further, to achieve
a similar result even in normed spaces of infinite dimension.
Provided we have a result of the Rivlin-Shapiro type, this is
possible for the approximation lemma extends to infinite
dimensional spaces, too.

Analogously to the above first proof, one may want to represent
an arbitrary convex body $K\subset X$ as intersection of balls;
but that does not always go through. The property of a normed
space $X$ that all convex bodies are intersections of closed
balls is called the Mazur Intersection Property, and this
property fails in some spaces: see e.g. \cite{Beer, GGS, GMP}.
Nevertheless, an even nicer proof of Lemma \ref{l:strictapprox}
can be presented if we apply the general fact, also well-known to
geometers, that strict convexity, in fact, is the dual property
to smoothness\footnote{We thank K{\'a}roly B{\"o}r{\"o}czky Jr.
reminding us to this idea, quite relevant in the present
context.} (cf. \cite{P}, e.g.), which gives way for the proof of
the general version.

{\it Second proof of Lemma \ref{l:strictapprox}.}\, Consider the
{\it dual convex body} $K^*\subset X^*$ of $K$, approximate it
closely by a {\it smooth} convex body $Q^*\subset X^*$, and then
take the pre-dual $Q\subset X$ of $Q^*$, which is then a strictly
convex body in $X$ and comes arbitrarily close to $K$ as $Q^*$
approximates $K^*$. To make everything explicit, one may argue by
adding a small ball to $K^*$, that is consider say $Q^*:=\delta
B^*+K^*$, and then let $\delta\to 0$ while taking the predual $Q$
of $Q^*$.

But to start with (to get a Rivlin-Shapiro type result, say), do
we also have a parallel supporting hyperplane lemma? Standard
proofs use compactness, which is no longer available in $X$.
Interestingly, even Lemma \ref{l:parsup} continues to hold in some
Banach spaces, in particular in reflexive Banach spaces, see
\cite[Proposition 2]{YRS}. That gives a way to recover the finite
dimensional, geometric proofs even in these normed spaces.
However, examples show\footnote{As written in \cite{YRS}, this
example was constructed by P. Wojtaszczyk.} (see \cite[Example
2]{YRS}), that the same assertion {\it fails} in some other Banach
spaces. Hence to settle the general case in a satisfactory way one
should combine our knowledge of $\alpha(K,x)$ more effectively. We
have the following result.

\begin{thm}[{\bf R{\'e}v{\'e}sz--Sarantopoulos, 2001, \cite{RS}}]
\label{th:RS} If $K\sbt X$ is an arbitrary convex body and if $ x
\in X\setminus K$ is arbitrary, then we have \eqref{CnTn}.
Moreover, $C_n(K, x )$ is actually a maximum, attained by
\eqref{PTnt}, where the notation is as in Theorem
\ref{th:RivlinShapiro}.
\end{thm}

Here we can observe that our results on the linear speed of growth
of $\al(K, x )$ (see \cite[Corollary 5.8., (5.17)]{RS}), together
with Theorem \ref{th:RS} give strong indications even for the
other Chebyshev problem, as we know that the Chebyshev polynomial
itself has leading coefficient $2^{n-1}$. Indeed, we have the
following result.

\begin{thm}[{\bf R{\'e}v{\'e}sz--Sarantopoulos, 2001 \cite{RS}}]\label{}
Let $K\sbt X$ be an arbitrary convex body and let $ v\in X$. Then
we have $$ A_n(K, v)={\frac{2^{2n-1}}{\tau{(K, v)}^n}}~,$$ and the
supremum is actually a maximum attained by a polynomial of the
form \eqref{PTnt} with some appropriately chosen $ v^*\in S^*$.
\end{thm}

Based on the determination of these extremal quantities, other
related questions were already addressed in approximation theory,
such as the uniqueness of the extremal polynomials, or the
existence of the so-called universal majorant polynomials. These,
in turn, have consequences e.g.\ concerning the approximation of
convex bodies by convex hulls of algebraic surfaces. For further
details we refer to \cite{Kr1} and \cite{R1}.

\section{Bernstein's Inequality}

If a univariate algebraic polynomial $p$ is given with degree
at most $n$, then by the classical Bernstein-Szeg\H o inequality
(\cite{SZ}, \cite{CS}, \cite{BO}) we have

\begin{equation}\label{classicbernstein}
|p'(x)| \leq \frac {n \sqrt{ ||p||_{C[a,b]}^2 - p^{2}({x}) }}
{\sqrt{(b-x)(x-a)}} \qquad (a<x<b).
\end{equation}

This inequality is sharp for every $n$ and every point $x\in
(a,b)$, as
$$
\sup \left\{ \frac {|p'(x)|}{\sqrt{ ||p||_{C[a,b]}^2 -
p^{2}({x})}} \, :\,\, \deg p \le n,\,|p(x)|< \|p\|_{C[a,b]}
\right\}= \frac {n}{\sqrt{(b-x)(x-a)}}\,\,.
$$

We may say that the upper estimate \eqref{classicbernstein} is
exact, and the right hand side is just the "true Bernstein
factor" of the problem.

In the multivariate setting a number of extensions were proved
for this classical result. However, due to the geometric variety
of possible convex sets replacing intervals of $\RR$, our present
knowledge is still not final. The exact Bernstein inequality is
known only for symmetric convex bodies, and we are within a bound
of some constant factor in the general, nonsymmetric case.

For more precise notation we may define formally for any
topological vector space $X$, a subset $K\subset X$, and a point
$x\in K$ the $n^{\rm th}$ "Bernstein factor" as

\begin{equation}\label{bernsteinfactor}
B_n(K,\bfx):=\frac 1n \sup \left\{ \frac {\|D p(\bfx)\|}{\sqrt{
||p||_{C(K)}^2 - p^{2}({\bfx})}} \, :\,\, \deg p \le
n,\,|p(\bfx)|<||p||_{C(K)} \right\}\,,
\end{equation}
where $Dp(\bfx)$ is the derivative of $p$ at $\bfx$, and even for
an arbitrary unit vector $\bfy\in X$
\begin{equation}\label{directionalbernsteinfactor}
B_n(K,\bfx,\bfy):=\frac 1n \sup \left\{ \frac {\langle D
p(\bfx),\bfy \rangle }{\sqrt{ ||p||_{C(K)}^2 - p^{2}({\bfx})}} \,
:\,\, \deg p \le n,\,|p(\bfx)|<||p||_{C(K)} \right\}\,.
\end{equation}

The perhaps nicest available method -- and, anyway, our favorite
-- is the {\em method of inscribed ellipses}, introduced into the
subject by Y. Sarantopoulos \cite{Sar}. This works for arbitrary
interior points of any, possibly nonsymmetric convex body.
However, other methods are in use and there is a striking
connection, only recently revealed, and still not fully
understood, between the method of inscribed ellipses and the
general approach through {\em pluripotential theory}. In this
survey we explain the method of inscribed ellipses, list the known
results, consider an instructive analysis of the case of the
simplex, and shortly comment on the intriguing questions still
open.

\subsection{}\label{sec:inscrellipse} Although for the reader's
convenience we include some short proofs, let us emphasize that,
unless otherwise stated, results in this section are due to
Sarantopoulos \cite{Sar}. The key of all of the method is the next

\begin{oldlemma}\label{inellipse}{\bf (Inscribed Ellipse Lemma,
Sarantopoulos, 1991).} Let $K$ be any subset in a vector space
$X$. Suppose that ${ x} \in K$ and the ellipse
\begin{equation}\label{ellipse}
{\bf r}(t) = \cos{t}~ \bfa + b \sin{t}~ { y} + {\bfx-\bfa} \qquad
(t \in [-\pi,\pi))\,.
\end{equation}
lies inside $K$. Then we have for any polynomial $p$ of degree at
most $n$ the Bernstein type inequality
\begin{equation}\label{Bernsteinellipse}
| \langle {D} p({ x}), { y} \rangle| \leq \frac{n}{b} \sqrt{
||p||_{C(K)}^2 - p^{2}({ x})}.
\end{equation}
\end{oldlemma}

{\it Proof.} Consider the trigonometric polynomial $T(t):= p({\bf
r}(t))$ of degree at most $n$. Since ${\bf r}(t) \subset K$ we
clearly have $ ||T|| \leq ||p||_{C(K)}$. According to the
Bernstein-Szeg\H o inequality \cite{SZ} (see also \cite{CS}) for
trigonometric polynomials,
$$
|T'(t)| \leq n \sqrt{ ||T||^2 - T(t)^2 }\leq n \sqrt{
||p||_{C(K)}^2 - p({\bf r}(t))^2 } \qquad (\forall t \in \RR).
$$
In particular, for $t = 0$, we get
$$
|T'(0)| \leq n \sqrt{||p_n||_{C(K)}^2 - p_n^{2}({\bf x}) }.
$$
By the chain rule
$$
T'(0) = \langle {D} p_{n}({\bf r}(0)), {\bf r}'(0) \rangle =
 \langle {D} p_{n}({  x}), b {  y} \rangle~,
$$
which completes the proof.

\begin{oldlemma}\label{unitball}{(\bf Sarantopoulos, 1991).}
Let $K$ be a centrally symmetric convex body in a vector space $X$
and ${ x} \in K$. The ellipse $ {\bf r}(t) = \cos{t}~ { x} + b
\sin{t}~ { y}~~ (t\in[-\pi,\pi))$ lies in $K$ whenever $$||{
y}||_{(K)} = 1~{\rm and}~ b = \sqrt{ 1 - ||{ x}||^{2}_{(K)} }. $$
\end{oldlemma}

{\it Proof.} The assertion is equivalent to $||{\bf r}(t)||_{(K)}
\leq 1$ for every $t$. By the triangle and Cauchy inequalities $$
||{\bf r}(t)||_{(K)}\! \leq\! |\cos{t}|~ ||{ x}||_{(K)} + b
|\sin{t}|~ ||{ y}||_{(K)}\! \leq \!\sqrt{\cos^2{t}+\sin^2{t} }
\sqrt{ ||{ x}||^{2}_{(K)} + b^2 ||{ y}||^{2}_{(K)}} \! =\! 1. $$
Lemma \ref{unitball} is proved.

Mutatis mutandis to the previous lemma we can deduce also the
following variant.

\begin{oldlemma}\label{symconvbody}{(\bf Sarantopoulos, 1991).}
Let $K$ be a centrally symmetric convex body in $X$, where $(X,
||\cdot||)$ is a normed space. Let $\varphi_{K}=\|\cdot\|_{(K)}$
be the Minkowski functional (norm) generated by $K$. Then for
every nonzero vector ${ y} \in X$ the ellipse $ {\bf r}(t) =
\cos{t}~ { x} + b \sin{t}~ { y}~~  (t\in[-\pi,\pi))$ lies in $K$
with $$ b := \frac{ \sqrt{1-\varphi^2(K, x)} }{ \varphi(K, y)}. $$
\end{oldlemma}

\begin{oldresult}\label{unball}{(\bf Sarantopoulos, 1991).}
Let $p$ be any polynomial of degree at most $n$ over the normed
space $X$. Then we have for any unit vector $ y \in X$ the
Bernstein type inequality
\begin{equation}\label{Bernsteineq}
| \langle {D} p({ x}), { y} \rangle| \leq \frac{n \sqrt{
||p||_{C(K)}^2 - p^{2}({ x}) }} {\sqrt{1-\| x\|^{2}_{(K)}}}.
\end{equation}
\end{oldresult}

{\it Proof.} The proof follows from combining Lemmas
\ref{inellipse} and \ref{unitball}.

\begin{oldresult}\label{scb}{(\bf Sarantopoulos, 1991).}
Let $K$ be a symmetric convex body and $ y$ a unit vector in the
normed space $X$. Let $p_n$ be any polynomial of degree at most
$n$. We have $$ | \langle {D} p_{n}({ x}), { y} \rangle| \leq
\frac{2n \sqrt{ ||p_n||_{C(K)}^2 - p_n^{2}({ x}) }} {\tau(K,{ y})
\sqrt{1 - \varphi^{2}(K, { x})}}. $$ In particular, with $w(K)$
standing for the width of $K$ we have $$ \|D p_n( x)\|\le \frac{2
n \sqrt{ ||p_n||_{C(K)}^2 - p_n^{2}({ x}) }} {w(K)
\sqrt{1-\varphi^{2}(K, { x})}}. $$
\end{oldresult}

{\it Proof.} Here we need to combine Lemmas \ref{inellipse} and
\ref{symconvbody} to obtain Theorem \ref{scb}.

\bigskip

It can be rather difficult to determine, or even to estimate the
$b$-parameter of the "best ellipse", what can be inscribed into a
convex body $K$ through $ x\in K$ and tangential to direction of $
y$. Still, we can formalize what we are after.
\begin{defn}\label{generalbestellipse}
For arbitrary $K\subset X$ and $ x\in K$, $ y\in X$ the
corresponding "best ellipse constants" are the extremal quantities
\begin{equation}\label{EKxy}
E(K, x, y):=\sup \{ b \,: \, \bfr \subset K \,\, {\hbox {\rm
with}}\,\, \bfr \,\, {\hbox {\rm as given in \eqref{ellipse}}} \,
\}
\end{equation}
and
\begin{equation}\label{EKx}
E(K, x):=\inf \{ E(K, x, y) \,: \,  y \in X , || y||=1\, \}\, .
\end{equation}
\end{defn}

Clearly, the inscribed ellipse method yields Bernstein type
estimates whenever we can derive some estimate of the ellipse
constants. In case of symmetric convex bodies, Sarantopoulos's
Theorems \ref{unball} and \ref{scb} are sharp; for the
nonsymmetric case we know only the following result.

\begin{oldresult}\label{KrooRevesz}
{\bf (Kro{\'o}--R{\'e}v{\'e}sz, \cite{KR}).} Let $K$ be an arbitrary
convex body, ${ x} \in {\rm int}K $ and $\|{ y}\| = 1$, where $X$
can be an arbitrary normed space. Then we have
\begin{equation}\label{krry}
| \langle {D}p( x),  y \rangle | \leq
 \frac {2 n \sqrt{ ||p||_{C(K)}^2 - p^{2}({ x}) }}
 { \tau(K, { y}) \sqrt{1 - \alpha(K,{ x})} } ,
 \end{equation}
for any polynomial $p$ of degree at most $n$. Moreover, we also
have
\begin{equation}\label{krrgrad}
||{D}~p({ x})|| \leq \frac {2n \sqrt{ ||p||_{C(K)}^2 - p^{2}({ x})
}} {w(K) \sqrt{ 1- \alpha(K, { x})} } \leq \frac {2 \sqrt 2 n
\sqrt{ ||p||_{C(K)}^2 - p^{2}({ x}) }} {w(K) \sqrt{ 1- \alpha^2(K,
{ x})} }.
\end{equation}
\end{oldresult}

Note that in \cite{KR} the best ellipse is not found; the
construction there gives only a good estimate, but not an exact
value of \eqref{EKxy} or \eqref{EKx}. In fact, here we quoted
\cite{KR} in a strengthened form: the original paper contains a
somewhat weaker formulation only.

One of the most intriguing questions of the topic is the following
conjecture, formulated first in \cite{RS}.

\begin{oldconjecture}\label{alphasquare}{\bf
(R{\'e}v{\'e}sz--Sarantopoulos).} Let $X$ be a topological vector
space, and $K$ be a convex body in $X$. For every point $ x \in
{\rm int} K$ and every (bounded) polynomial $p$ of degree at most
$n$ over $X$ we have $$ \|D p( x)\|\le \frac{2 n \sqrt{
||p||_{C(K)}^2 - p^{2}({ x}) }} {w(K) \sqrt{1-\alpha^{2}(K,
x)}}\,, $$ where $w(K)$ stands for the width of $K$.
\end{oldconjecture}

\subsection{}\label{sec:simplex} We denote $|x|_2 :=
(\sum_{i=1}^{d} x_{i}^2)^{1/2}$ the Euclidean norm of $x = (x_{1},
\ldots, x_{d}) \in {\mathbb R}^{d}$. Let
$$
\Delta := \Delta_d := \{(x_{1},\ldots,x_{d}) : x_{i} \geq 0,
i=1,\ldots,d, \sum_{i=1}^{d}x_{i} \leq 1\}
$$
be the standard simplex in ${\mathbb R}^{d}$. For fixed $ { x} \in
{\rm int} \Delta, $ and $ { y} = (y_{1},\ldots,y_{d}),~ |{ y}|_{2}
= 1$ the best ellipse constant of $\Delta$ is, by Definition
\ref{generalbestellipse}, $ E(\Delta, x, y) $. By a tedious
calculation via the Kuhn-Tucker theorem and some geometry, the
following was obtained in \cite{MR}.
\begin{prop}[{\bf Milev-R{\'e}v{\'e}sz, 2003}]\label{bestb}
We have
\begin{equation}\label{Evalue}
E(\Delta, x, y) = \left\{ \frac{y_{1}^2}{x_{1}} + \cdots +
\frac{y_{d}^2}{x_{d}} + \frac{(y_{1} + \ldots + y_{d})^2}{1 -
x_{1} - \ldots - x_{d}} \right\}^{-1/2}.
\end{equation}
\end{prop}

\begin{thm}[\bf Milev-R{\'e}v{\'e}sz, 2003]\label{ellipseyield} Let $p_{n} \in {\mathcal P}_{n}^{d}$.
Then for every ${  x} \in {\rm int} \Delta$ and ${  y} \in
{\mathbb S}^{d-1}$ we have
\begin{equation}\label{directionalyield}
| D_{{  y}} p_{n}({  x}) | \leq \frac {n \sqrt{
||p_n||_{C(\Delta)}^2 - p_n^{2}({  x}) }} { E(\Delta,{  x}, { y})
},
\end{equation}
where $E(\Delta,{  x}, {  y})$ is as given in \eqref{Evalue}.
\end{thm}

>From now on let us restrict ourselves to the case $d=2$. We
denote the vertices of $\Delta$ by $O = (0,0), A = (1,0), B =
(0,1)$ and the centroid (i.e. mass point) of $\Delta$ by $M =
(1/3,1/3)$. A calculation shows that $1 - \alpha(\Delta, {  x}) =
2 r( x)$, with
\begin{equation}\label{r}
r := r( x) := \min \{ x_{1}, x_{2}, 1 - x_{1} - x_{2}\}
               = \left\{
           \begin{array}{ll}
           x_{1}, & {  x} \in \triangle OMB \\
            x_{2}, & {  x} \in \triangle OMA \\
           1 - x_{1} - x_{2}, & {  x} \in \triangle AMB \\
           \end{array} \right.
\end{equation}
and if ${  y} = (\cos \varphi, \sin \varphi)$ $(0 \leq \varphi
\leq \pi)$ then $$ \tau(\Delta, {  y}) = \left\{
      \begin{array}{ll}
      1/(y_{1} + y_{2}),  &  \varphi \in [0,\pi/2] \\
      1/y_{2}, & \varphi \in ( \pi/2, 3\pi/4] \\
      -1/y_{1}, & \varphi \in ( 3\pi/4, \pi]. \\
      \end{array}
      \right.
$$

Note that the inequality
\begin{equation}\label{oldright7}
{  \frac{1} { E(\Delta,{  x}, {  y})} }    \leq { \frac {2}
{\tau(\Delta, {  y}) \sqrt{1 - \alpha(\Delta, {  x})} } }
\end{equation}
holds true for every $x \in {\rm int} \Delta$ and $y \in {\mathbb
S}^1$, i.e. estimate \eqref{directionalyield} is better than
\eqref{krry} when $K = \Delta$. Accordingly, we can derive a new
estimation for $Dp_n(x)$.

\begin{prop}[\bf Milev-R{\'e}v{\'e}sz, 2003]\label{gradient}
Let $p_{n} \in {\mathcal P}_{n}^{2}$. Then for every ${  x} \in
{\rm int} \Delta$ we have
\begin{equation}\label{oldright11}
 | {D}p_{n}({  x}) |_{2} \leq n E({  x})\sqrt{
||p_n||_{C(\Delta)}^2 - p_n^{2}({  x}) },
\end{equation}
where
\begin{equation}\label{oldright12}
 E({  x}) =
\sqrt{ \frac{x_{1}(1-x_{1}) + x_{2}(1-x_{2}) +  D({  x}) }
   { 2 x_{1} x_{2} (1 - x_{1} - x_{2}) }}
\end{equation}
with $$ D({  x}) =  \sqrt{ [x_{1}(1-x_{1}) + x_{2}(1-x_{2})]^2 - 4
x_{1} x_{2} (1- x_{1} - x_{2})}. $$
\end{prop}
Note that the inequality $$ [x_{1}(1-x_{1}) + x_{2}(1-x_{2})]^2 -
4 x_{1} x_{2} (1- x_{1} - x_{2}) > [x_{1}(1-x_{1}) -
x_{2}(1-x_{2})]^2 $$ holds true for ${  x} \in {\rm int} \Delta$,
hence $D({  x}) > 0$.

Using this estimate, an improvement of the constant $2$ to
$\sqrt{3}$ was achieved in Theorem \ref{KrooRevesz} for the
special case of $K = \Delta$, c.f. \cite{MR}. Of more interest is
the next estimate comparing to the conjectured quantity with $1 -
\alpha^2(\Delta, x)$.
\begin{thm}[\bf Milev-R{\'e}v{\'e}sz, 2003]\label{squarecompare}
Let $p_{n} \in {\mathcal{P}}_{n}^{2}$ and
 $||p_{n}||_{C(\Delta)} =
1$. Then for every ${  x} \in {\rm int}~\Delta$ we have
\begin{equation}\label{number}
 |{D}p_{n}({  x})|_{2} \leq
 \frac{ \sqrt{3+\sqrt{5}}~ n
 \sqrt{||p_n||_{C(\Delta)}^2 - p_n^{2}({  x})} }
       { w(\Delta) \sqrt{ 1 - \alpha^2(\Delta, {  x})}}.
\end{equation}
\end{thm}
It was checked that this is the most what follows from the
inscribed ellipse method, interpreted as considering
$E(\Delta,x)$ the exact yield it gives.

This improves the constant in Theorem \ref{KrooRevesz} but falls
short of Conjecture \ref{alphasquare}, since $2\sqrt{2} =
2.8284\ldots > \sqrt{3+\sqrt{5}} = 2.2882\ldots > 2.$

\subsection{}\label{sec:ridge}
Let us consider the following question. All known lower estimates
for the Bernstein factors used some kind of ridge polynomials,
i.e. polynomials composed from a linear form and some (in fact, a
Chebyshev) polynomial. Can one sharpen these lower estimates to
the extent that Conjecture \ref{alphasquare} will be disproved?

Recall that ridge polynomials are defined as $$ {\mathcal
R}_n:=\{p\in {\mathcal P} \,:\, p( x)=P(L( x)), L\in X^{*},
P\in{\mathcal P}_n(\RR) \}, \qquad {\mathcal R}:=
\bigcup_{n=1}^{\infty}{\mathcal R}_n. $$ By easy linear
substitution we may assume that ridge polynomials are expressed by
using some $L(x)=t(K, v^{*}, x)$, as defined in \eqref{lineart}.

\begin{defn}\label{C(x,y)}
For any $n\in\NN$ the corresponding "ridge Bernstein constant" is
$$ C_{n}(K,{  x}, {  y}) := \frac{1}{n} \sup_{ R \in  {\mathcal
R}_n, |R({  x})|<\|R\|_{C(K)} } \frac{|\langle DR( x),  y
\rangle|} {\sqrt{\|R\|_{C(K)}^2-R^2({  x})}}\,\,. $$
\end{defn}

\begin{prop}[\bf Milev-R{\'e}v{\'e}sz, 2003]\label{chebyshevcomp}
For every convex body $K$ and ${  x} \in {\rm int} K$, ${  y} \in
S^{*}$ we have $$ C_{n}(K, {  x}, {  y}) \leq {\frac 2 {\tau(K, {
y})}} {\frac 1 {\sqrt{1 - \alpha^2(K, {  x})}}}\,\,. $$
\end{prop}

{\it Proof.} By the chain rule we have for any $R\in {\mathcal
R}_n, \,\, R=P (t( x))$ the formula
$$
| \langle DR({  x}), { y} \rangle |
 = \left| P'(t({  x})) {\frac 2{w(K,{  v}^*)}}
 \langle {  v}^*, {  y} \rangle \right|
 \leq {\frac 2{\tau(K, {  y})}} \left| P'(t({  x})) \right|.
$$
Applying the Bernstein-Szeg\H o inequality for $s \in (-1,1)$
we get
$$
\frac{|P'(s)|} {\sqrt{\|P\|_{C[-1,1]}^2-P^2(s)}} \le
{\frac n {\sqrt{1 - s^2}}}\,\,\, .
$$
Note that for $T_n$, the classical Chebyshev polynomial of degree
$n$, (and only for that) this last inequality is sharp. Putting
$s:=t( x)$ and combining the previous two inequalities we are led
to
$$
 {\frac1n} \frac{|\langle D R( x),  y \rangle |}
 {\sqrt{\|R\|_{C(K)}^2-R({  x})}} \leq \frac{2}{\tau(K,{  y})}
 \frac{1}{\sqrt{1 -  t^{2}(K,{v}^{*}, {  x})}}\,\,.
$$
Taking supremum with respect to ${  v}^{*}\in S^{*}$ on the right
hand side, we obtain a bound independent of ${  v}^{*}$. In fact,
according to \eqref{lambdadef} and \eqref{alphalambda} (see also
\cite[Proposition 4.1]{RS},), the supremum is a maximum and is
equal to $ \frac 2{\tau(K,{  y})} \frac 1 {\sqrt{1 - \alpha^2(K,
{  x})}}\,\,$. Thus taking supremum also on the left hand side,
Theorem \ref{chebyshevcomp} obtains.

\medskip
Whence ridge polynomials satisfy Conjecture \ref{alphasquare},
always.

\medskip
It follows from the definitions and Lemma \ref{inellipse} that
$$
C_{n}(K,{  x}, {  y}) \leq B_{n}(K,{  x}, {  y}) \leq \frac{1}
{E(K,{  x}, {  y})}\,\,\,.
$$

For the case of the standard simplex we have a converse
inequality.
\begin{prop}[\bf Milev-R{\'e}v{\'e}sz, 2003]\label{converseineq}
For every ${  x} \in {\rm int} \Delta$ and ${  y} \in {\mathbb
S}^{d-1}$ we have the inequality $$ \frac{1} {E(\Delta,{  x}, {
y})} \leq \sqrt{d}~
 C_{n}(\Delta, {  x}, {  y})\,\,.
$$
\end{prop}

\begin{cor}[\bf Milev-R{\'e}v{\'e}sz, 2003]\label{bncn}
For every ${  x} \in {\rm int} \Delta$ and ${  y} \in {\mathbb
S}^{d-1}$ we have $$
 1 \leq \frac{B_{n}(\Delta,{  x}, {  y})}
{C_{n}(\Delta, {  x}, {  y})} \leq \sqrt{d}, \qquad 1 \leq
\frac{B_{n}(\Delta,{  x})} {C_{n}(\Delta, {  x})} \leq
\sqrt{d}\,\,. $$
\end{cor}

Note that in the paper \cite{KR} it is proved that for every $
x_{0} \in {\rm int}K$ there is a direction $ y_{0}$ and a ridge
polynomial $T_n(K, v^*_0, x)$ such that
$$
\frac{1}{n} \frac{D_{y_0}T_n(K, v^*_0, x_0)} {\sqrt{1-T_n(K,
v^*_0, x_0)^2}} = \frac{2} { \tau(K, {  y}_{0}) \sqrt{1 -
\alpha^{2}(K,{ x}_{0})} }\,\,.
$$
Consequently,
$$
C_{n}(K,
 x_{0},  y_{0}) \geq \frac{2} { \tau(K, {  y}_{0}) \sqrt{1 -
\alpha^{2}(K,{  x}_{0})} }\,\,.
$$
Hence, for every $ x_0\in {\rm int}K$ there is a $y_0$ such that
$$
\frac{B_n (K, x_0, y_0)} {C_{n}(K, x_0, y_0)} \leq \sqrt{2}\,.
$$
Comparing this to Corollary \ref{bncn}, we see that (for the case
of the simplex) the latter ratio remains uniformly bounded for
all $ x$ and $ y$.

\subsection{}\label{sec:pluripot} Another method of considerable
success in proving Bernstein (and Markov) type inequalities is the
pluripotential theoretical approach. Classically, all that was
considered only in the finite dimensional case, but nowadays even
the normed spaces setting is cultivated. To explain these, one
needs an understanding of {\em complexifications of real normed
spaces}, see e.g. \cite{MST, BARAN4}, as well as the {\em
Siciak-Zaharjuta extremal function} $V(z)$. In fact, the latter,
by the celebrated Siciak-Zaharjuta Theorem, can be expressed both
by {\em plurisubharmonic functions} from the {\em Lelong class},
and also just by logarithms of the absolute values of polynomials.
We spare the reader from the first, referring to \cite{Kli} as a
general, nice introduction to pluripotential theory, and restrict
ourselves to the latter, perhaps easier to digest formulation.
That is very much like the Chebyshev problem \eqref{Chebyfindef}
in \S \ref{sec:Cheby}, except that we consider it all over the
complexification $Y:=X+iX$ of $X$, take logarithms, and after
normalization by the degree, merge the information derived by all
polynomials of any degree into one clustered quantity. Namely, for
any bounded $E\subset Y$ $V_E$ vanishes on $E$, while outside $E$
we have the definition
\begin{equation}\label{SZV}
V_E(z):=\sup \{\frac{1}{n}\log|p(z)|~:~ 0\neq p \in \PPn(Y),~
||p||_E\le 1,~n\in \NN \} \qquad (z\notin E)
\end{equation}
For $E\subset X$ one can easily restrict even to $p\in\PP(X)$. For
the theory related to this function and some recent developments
concerning Bernstein and also Markov type inequalities for convex
bodies or even more general sets, we refer to \cite{BARAN1,
BARAN2, BARAN3, BARAN4, Kli, MAU, LEV, PP, Ple}.

Now consider $E=K\subset X$, where $K$ is now a convex body. Our
more precise result in Theorem \ref{th:RS} yields
$V_K(x)=\sup_{n\in\NN} \log|C_n(K,x)|/n = \lim_{n\in\NN}
\log|C_n(K,x)|/n = \alpha(K,x)+\sqrt{\alpha(K,x)^2-1}$, as an easy
calculation with the last expression in \eqref{Tndef} shows
together with the fact that $\log|C_n(K,x)|/n$ increases. However,
in the Bernstein problem the values of $V_K$ are much more of
interest for complex points $z=x+iy$, in particular for $x\in K$
and $y$ small and nonzero. More precisely, the important quantity
is the normal (sub)derivative
\begin{equation}\label{normalder}
D_{y}^{+} V_E(x):=\liminf_{\epsilon\to 0}\, \frac{V_E(x+i\epsilon
y)}{\epsilon}~,
\end{equation}
as this quantity occurs in the next estimation of the directional
derivative.
\begin{thm}[\bf Baran, 1994 \& 2004]\label{th:BaranBernstein}
Let $E\subset X$ be a bounded, closed set, $x\in \intt E$ and
$0\ne y\in X$. Then for all $p\in \PPn(X)$ we have
\begin{equation}\label{Baranestimate}
|\langle Dp(x), y \rangle|\le n D_{y}^{+} V_E(x)
\sqrt{||p||^2_E-p(x)^2}~.
\end{equation}
\end{thm}
In fact, \cite{BARAN1} contains this only for $\RRd$ and partial
derivatives, but by applying rotations of $E$, all directional
derivatives follow; the case of infinite dimensional spaces are
considered in \cite{BARAN4}. See also \cite{BLW, R}.

It is not obvious, how such estimates can be applied to concrete
cases. First, one has to find the precise value of $V_E$, in such
a precision, that even the derivative can be computed: then the
derivatives must be obtained and only then do we really have
something. However, even that is addressed by considering the
Bedford-Taylor theory of the Monge-Ampere equation and the
equilibrium measure \cite{BT}, as the density of the equilibrium
measure gives the extremal function. In some concrete applications
all that may be calculated. A particular example (see
\cite{BARAN0}, \cite[Example 5.4.7]{Kli}, \cite[Example
4.8]{BARAN3}) is the following.
\begin{prop}[\bf Baran, 1988]\label{Baranexample} The extremal
function of the standard simplex in $\RRd$ is $V_{\Delta}(z)=\log
h(|z_1|+|z_2|+\dots+|z_n|+|1-(z_1+z_2+\dots+z_n)|)$.
\end{prop}
>From this and the calculation with the rotated directions etc, we
calculated in \cite{R} the following surprising corollary.
\begin{prop}\label{Reveszexample} The
above pluripotential theoretical estimate of Baran gives for the
standard triangle of $\RR^2$ the result exactly identical to
\eqref{directionalyield}.
\end{prop}
Much remains to explain in this striking coincidence, the first
being the next.
\begin{hypothesis}\label{hyp:allinone} Let $K\subset X$ be a convex body.
Then for all points $x\in \intt K$ the inscribed ellipse method
and the pluripotential theoretical method of Baran results in
exactly the same estimate, i.e. for all $y\in S^*$ we have
\begin{equation}\label{ellipseeqppot}
D_{y}^{+} V_K(x) = \frac {1}{E(K,x,y)}~.
\end{equation}
\end{hypothesis}

All people like to believe that his method(s) are the ultimate
ones. However, it is quite unclear which one is the right one in
the Bernstein problem. If any of the inscribed ellipse method or
the pluripotential theoretical method of Baran is right - or, in
case of validity of Hypothesis \ref{hyp:allinone}, if both are
precise - then Conjecture \ref{alphasquare} would fail. Still, it
seems worthy to formulate these contradictory assumptions.

\begin{hypothesis}\label{hyp:Baranright} Let $K\subset X$ be convex body.
Then for all points $x\in \intt K$ the exact Bernstein factor is
just what results from the pluripotential theoretical method of
Baran:
\begin{equation}\label{pluriexact}
B_n(K,x) = \sup_{y\in S^{*}} D_{y}^{+} V_K(x)~.
\end{equation}
\end{hypothesis}

\begin{hypothesis}\label{hyp:ellipseright} Let $K\subset X$ be convex body.
Then for all points $x\in \intt K$ the exact Bernstein factor is
just what results from  the inscribed ellipse method of
Sarantopoulos:
\begin{equation}\label{ellipseexact}
B_n(K,x) = \frac {1}{E(K,x)}~.
\end{equation}
\end{hypothesis}

Note that we already know that these hypothesis are certainly not
true for the directional derivatives of {\em all} directions $y\in
S^*$, where both methods can be improved upon for some $y$, see
\cite{R0}. Care has to be exercised in formulating conjectures and
hypothesis in these matters: the situation is more complex than
one might like to have, and the simple heuristics of extending the
results of the symmetric case do fail sometimes. In this respect
see also \cite{BLM, BLW, LEV, MAU}

For some other interesting assertions and conjectures, (sometimes
more addressed to the pluripotential theoretical aspects than the
Bernstein inequality itself), and an analysis of them we refer to
\cite{BARAN4, R0}.

Also, another real, geometric method, of obtaining Bernstein type
inequalities, due to Skalyga \cite{S1, S3}, is to be mentioned
here: the difficulty with that is that to the best of our
knowledge, no one has ever been able to compute, neither for the
seemingly least complicated case of the standard triangle of
$\RR^2$, nor in any other particular non-symmetric case the yield
of that abstract method. Hence in spite of some remarks that the
method is sharp in some sense, it is unclear how close these
estimates are to the right answer and what use of them we can
obtain in any concrete cases.

\section{Further inequalities and problems for
solution}\label{sec:further}

\subsection{}\label{sec:Markov} Finally let us touch upon a few
other questions and problems generally in the center of interest
for approximation theorists. One is the so-called Markov problem,
which is the question of obtaining {\em uniform} estimates, (as
opposed to pointwise ones in the Bernstein problem), to the size
of the gradient vector of a polynomial all over the convex body
$K\subset X$. Note that while the Bernstein-Szeg\H o type
estimates are quite good for a given point $x$, their use is less
and less towards the boundary: in fact, at the boundary the
estimate tends to infinity. This is so even in the one dimensional
case of $\RR$, and is inherent in the problem, due to the
improvement, generally valid only inside the body, with respect to
dependence on the degree of the polynomial. Indeed, $B_n(x)$ was
normalized just by $n$, the degree, while the classical Andrei
Markov inequality
\begin{equation}\label{Markovfirst}
\|p'\|_{[a,b]} \le \frac{2n^2}{b-a} \|p\|_{[a,b]}
\end{equation}
is sharp, excluding a "uniform Bernstein inequality" even in
$\RR$.

For symmetric convex bodies, also by the above described method of
inscribed ellipses, Sarantopoulos was able to obtain that
\begin{equation}\label{Markovsymm}
\Vert Dp \Vert_K \le \frac{2n^2 \|p\|_K}{w(K)}~,
\end{equation}
or, in case $\|\cdot\|=\|\cdot\|_{(K)}$, i.e. when $K$ is the unit
ball of the normed space $X$, then
\begin{equation}\label{Markovball}
\Vert Dp \Vert \le n^2 \|p\|~,
\end{equation}
a fully satisfactory answer for the symmetric case. A nice,
elementary argument of Wilhelmsen \cite{Wilh} presented
\eqref{Markovsymm} in the full generality of convex, not
necessarily symmetric bodies, but with a factor 4 in place of 2.
It was shown in \cite{NS} that \eqref{Markovsymm} does not remain
true for all convex bodies. Finally, Skalyga \cite{S2, S3} found
\begin{equation}\label{Markovgeneral}
\Vert Dp \Vert_K \le \frac{2n \cot\left(\frac{\pi}{4n}\right)
\|p\|_K}{w(K)}~.
\end{equation}
Note that
\begin{equation}\label{Skalygaass}\notag
2n \cot\left(\frac{\pi}{4n}\right) ~ \sim \frac8{\pi} n^2 \qquad
\text{as}\quad n\to \infty~.
\end{equation}
The author of \cite{S2} also remarks that the estimate
\eqref{Markovgeneral} is sharp in the sense of being subject to no
improvement in the full generality of all convex bodies and in all
normed spaces. For some further information see \cite{A, BARAN4,
BCG, BO, Kr2, PP, Ple}.

\subsection{}\label{sec:higherMarkov} It is important, in
particular for doing analysis on infinite dimensional spaces, to
have a control over the size of derivatives of any order. In the
classical case of a real interval it was done by Vladimir Markov,
and the answer e.g. for $I=[-1,1]$ is
\begin{equation}\label{Markovhigh}
\|p^{(k)}\|_I \le T_n^{(k)}(1) \|p\|_{I} =
\frac{(n^2(n^2-1)\dots(n^2-(k-1)^2)}{(2k-1)!} \|p\|_{I}\qquad
(k=1,2,\dots,n) ~,
\end{equation}
which is sharp again for the Chebyshev polynomial $T_n$.

At present we are far from having a nearly as precise estimate as
\eqref{Markovhigh} for the general case of normed spaces. These
can not be obtained, not even for dimension one, by simple
iterations of the estimates for the first derivative, what gives
only substantially weaker results. However, for the important
special case of a Hilbert space an exact extension of this
inequality is known, see \cite{MS}.

\begin{oldresult}{\bf (Munoz--Sarantopoulos, 2002).}\label{th:Hilbert}
Let $H$ be any Hilbert space and $p$ be an arbitrary polynomial of
degree $n$, that is $p\in\PPn(H)$. Then we have
\begin{equation}\label{HilbertMarkov}
|\|p^{(k)}\| \le T_n^{(k)}(1) \|p\| \qquad (k=1,2,\dots,n) ~.
\end{equation}
\end{oldresult}

Harris has a number of results on Bernstein and Markov
inequalities related also to higher order derivatives, see
\cite{H2}, and the extremely readable survey \cite{H1} in
particular. One would like to decide the following.
\begin{oldconjecture}{\bf (Harris).}\label{conj:Harris}
Let $X$ be any Banach space and $p\in\PPn(X)$. Then
\eqref{HilbertMarkov} holds true.
\end{oldconjecture}
It seems that neither the inscribed ellipse, nor the
pluripotential theoretic methods above can be applied to higher
derivatives, at least not directly. Hence even in the symmetric
(i.e., norm unit ball) case there is no obvious way to get close
to the conjecture (or the truth).

\subsection{}\label{sec:Schur} As mentioned above, combining or
iterating known estimates does not necessarily give best results
even if the parts put together are exact in their kind. Another
example is the following classical question, which can be
considered a composition of the Bernstein problem and the
Chebyshev problem (although simple and basic in itself). The
question is that how large can the {\em derivative} of a
polynomial be at point $x$ not inside, but outside (may be
distant) of the set of normalization. The classical version for
dimension 1 was already known to Chebyshev and reads
\begin{equation}\label{Schurone}
\max \{p^{(k)}(x) ~:~ p\in\PPn(\RR),~ ||p||_{[-1,1]}\le
1\}=T^{(k)}(x)~,
\end{equation}
another extremal property of the classical Chebyshev polynomials,
see e.g. \cite[p. 93]{Ri}. This classical inequality can easily be
obtained by considering Lagrange interpolation of $p$ on the nodes
of maxima of the Chebyshev polynomials, that is at the point
system $\{\cos(k\pi/n)\}^n_{k=0}$. But what is the answer for the
similar question in the multivariate case? Since now directional
derivatives $D_y p(x)= \langle Dp(x),y \rangle$ of all directions
$y \in S$ occur, the problem does not reduce to a one dimensional
question. It would be interesting, but non-trivial, to settle this
question even in case $k=1$.

\end{document}